\documentclass{article}
\usepackage{graphicx}
\usepackage{graphics, graphicx, latexsym, amssymb}
\newtheorem{propo}{Proposition}
\newtheorem{coro}{Corollary}

\begin{document}

\title{Polyominoes with nearly convex columns: \\
A model with semidirected blocks}
\author{Svjetlan Fereti\'{c} \footnote{e-mail: svjetlan.feretic@gradri.hr} \\ 
Faculty of Civil Engineering, University of Rijeka, \\ 
Viktora Cara Emina 5, 51000 Rijeka, Croatia}
\maketitle

\begin{abstract}
In most of today's exactly solved classes of polyominoes, either all members are convex (in some way), or all members are directed, or both. If the class is neither convex nor directed, the exact solution uses to be elusive. This paper is focused on polyominoes with hexagonal cells. Concretely, we deal with polyominoes whose columns can have either one or two connected components. Those polyominoes (unlike the well-explored column-convex polyominoes) cannot be exactly enumerated by any of the now existing methods. It is therefore appropriate to introduce additional restrictions, thus obtaining solvable subclasses. In our recent paper, published in this same journal, the restrictions just mentioned were semidirectedness and an upper bound on the size of the gap within a column. In this paper, the semidirectedness requirement is made looser. The result is that now the exactly solved subclasses are larger and have greater growth constants. These new polyomino families also have the advantage of being invariant under the reflection about the vertical axis.
\end{abstract}

\vspace{8mm}

\vspace{3mm}

\vspace{3mm}
\noindent \textit{Suggested running head:} Polyominoes with nearly convex columns

\newpage

\section{Introduction}

In our previous paper \cite{semi}, we began to search for polyomino models which are more general than column-convex polyominoes, but still have reasonably simple \textit{area generating functions}. In \cite{semi}, we introduced \textit{level} $m$ \textit{cheesy polyominoes} ($m=1,\: 2,\: 3,\ldots$), and here we shall introduce another sequence of models, which we call \textit{level} $m$ \textit{polyominoes with cheesy blocks} ($m=1,\: 2,\: 3,\ldots$). 

At every level, both cheesy polyominoes and polyominoes with cheesy blocks have a rational area generating function. However, at any given level, cheesy polyominoes are a rather small subset of polyominoes with cheesy blocks. The latter set of polyominoes has a greater \textit{growth constant} than the former set. For example, the growth constant of level one cheesy polyominoes is $4.114907\ldots \:$, while the growth of level one polyominoes with cheesy blocks is $4.289698\ldots \:$. (By the growth constant we mean the limit $\lim_{n\rightarrow\infty} \sqrt[n]{a_n}$, where $a_n$ denotes the number of $n$-celled elements in a given set of polyominoes.)
In addition, if we reflect a polyomino with cheesy blocks about the vertical axis, we get a polyomino with cheesy blocks again. This kind of invariance under reflection is enjoyed by column-convex polyominoes, but not by cheesy polyominoes. Admittedly, counting level $m$ polyominoes with cheesy blocks requires some more effort than counting level $m$ cheesy polyominoes. Anyway, at level one, polyominoes with cheesy blocks are not very hard to count. In this paper, the level one model is solved in full detail, the solution of the level two model is outlined, and the result for the level three model is stated with no proof. (Just as with cheesy polyominoes, as level increases, the computations quickly gain in size.) Our computations are done by using Bousquet-M\'{e}lou's \cite{Bousquet} and Svrtan's~\cite{Svrtan} ``turbo" version of the Temperley method \cite{Temperley}. 

If the reader is interested in the history of polyomino enumeration, or in the role which polyominoes play in physics and chemistry, then he/she may refer to the recently published book \cite{book}. However, if the reader would settle for a few lines, then it could be enough to see the introduction of our previous paper \cite{semi}.

\section{Definitions and conventions}

Some of the relevant definitions were already stated in the section ``Definitions and conventions" of our previous paper \cite{semi}. Those ``old" definitions are not repeated here because, whether we repeat them or not, it is natural to read \cite{semi} before reading this paper.

In this paper, we deal with polyominoes with hexagonal cells. When we write ``a polyomino", we actually mean ``a hexagonal-celled polyomino". 

Suppose that $P$ is such a polyomino that the first (\textit{i.e.}, leftmost) column of $P$ has no gap and that, in every pair of adjacent columns of $P$, every connected component of the right column has at least one edge in common with the left column. Then we say that $P$ is a \textit{rightward-semidirected polyomino}.

A polyomino $P$ is a level $m$ cheesy polyomino if the following holds:

\begin{itemize}
\item $P$ is a rightward-semidirected polyomino,
\item every column of $P$ has at most two connected components,
\item if a column of $P$ has two connected components, then the gap between the components consist of at most $m$ cells.
\end{itemize}

See Figure 1.

\begin{figure}
\begin{center}
\includegraphics[width=55mm]{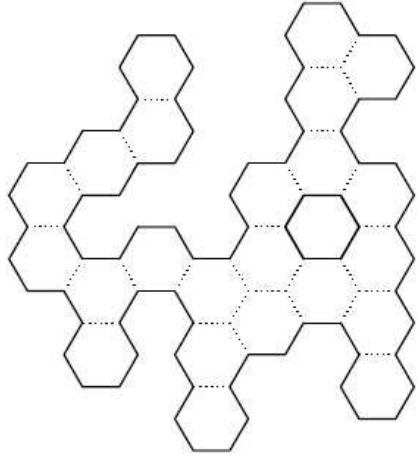}
\caption{A level one cheesy polyomino.}
\end{center}
\end{figure}

Suppose that $P$ is such a polyomino that the last (\textit{i.e.}, rightmost) column of $P$ has no gap and that, in every pair of adjacent columns of $P$, every connected component of the left column has at least one edge in common with the right column. Then we say that $P$ is a \textit{leftward-semidirected polyomino}.

A polyomino $P$ is a \textit{bird} if the following holds:

\begin{itemize}
\item $P$ has exactly one gap-free column (\textit{a}, say),
\item if we take the union of $a$ and of all (zero or more) columns lying to the left of $a$, the result is a leftward-semidirected polyomino,
\item if we take the union of $a$ and of all (zero or more) columns lying to the right of $a$, the result is a rightward-semidirected polyomino.
\end{itemize}

See Figure 2.

\begin{figure}
\begin{center}
\includegraphics[width=62.5mm]{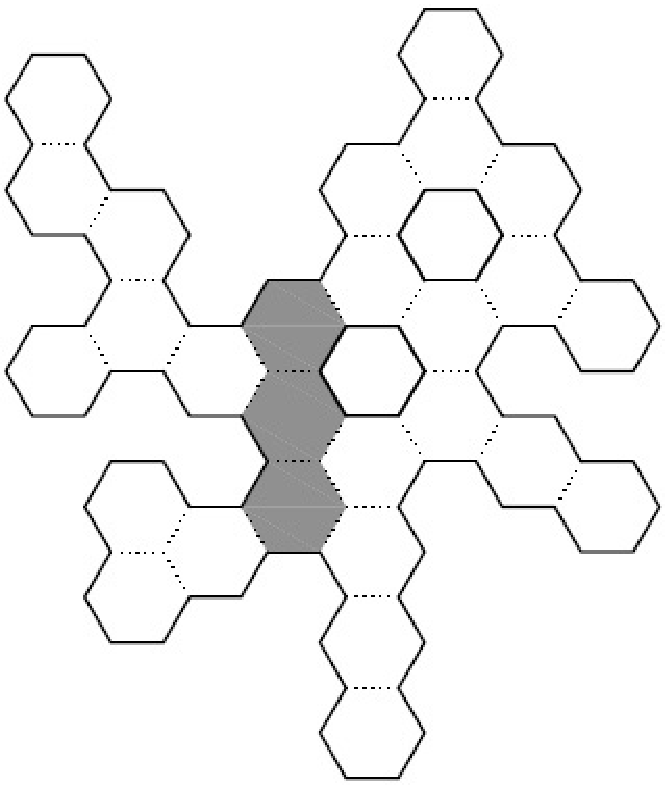}
\caption{A bird.}
\end{center}
\end{figure}

A polyomino $P$ is a \textit{level} $m$ \textit{polyomino with cheesy blocks} if the following holds:

\begin{itemize}
\item there exist positive integers $k$ and $i_1,\: i_2,\ldots ,\: i_k$ such that the first $i_1$ columns of $P$ form a bird, the next $i_2$ columns of $P$ form a bird,$\ldots,\: $ the last $i_k$ columns of $P$ form a bird,
\item every column of $P$ has at most two connected components,
\item if a column of $P$ has two connected components, then the gap between the components consists of at most $m$ cells.
\end{itemize}

See Figure 3.

\begin{figure}
\begin{center}
\includegraphics[width=85mm]{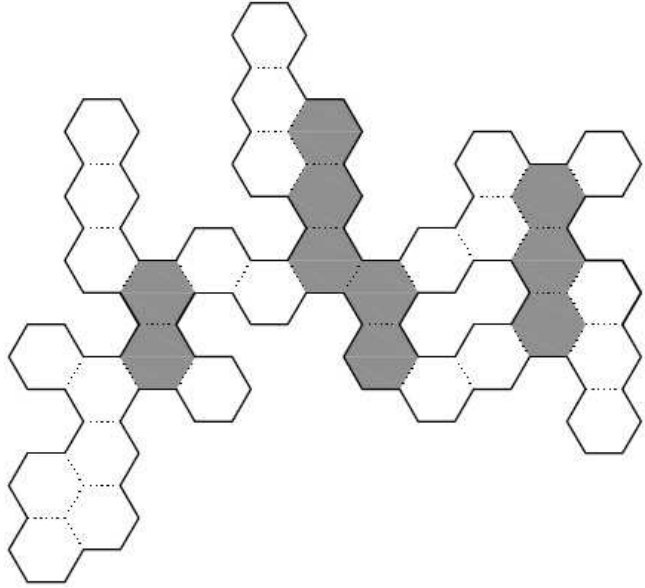}
\caption{A level one polyomino with cheesy blocks.}
\end{center}
\end{figure}

Notice that a polyomino with cheesy blocks may have more than one decomposition into birds. Nevertheless, this ambiguity will not bother us during the enumeration. Namely, our enumeration method is not a bird-by-bird one, but rather a column-by-column one. Next, it is easy to see that, if $P$ is a cheesy polyomino, then $P$ is also a polyomino with cheesy blocks. For example, say that a cheesy polyomino $P$ has $10$ columns and that the $1$st, $4$th and $5$th columns are gap-free, whereas each of the remaining $7$ columns does have a gap. Then the first $3$ columns form a bird, the $4$th column itself forms a bird, and the last $6$ columns form a bird as well.

Furthermore, if we reflect a bird about the vertical axis, we get a bird again. Therefore, if we reflect a polyomino with cheesy blocks about the vertical axis, we get a polyomino with cheesy blocks again.

If a polyomino $P$ is made up of $n$ cells, we say that the \textit{area} of $P$ is $n$.

Let $a$ be a column of a polyomino $P$. By the \textit{height} of $a$ we mean the number of those cells which make up $a$ plus the number of those (zero or more) cells which make up the gaps of $a$. For example, in Figure 2, the highlighted column has height $3$, and the next column to the left has height $3$, too. 

Let $R$ be a set of polyominoes. By the \textit{area generating function} of $R$ we mean the formal sum

\begin{displaymath}
\sum_{P \in R} q^{area\ of\ P}.
\end{displaymath}

By the \textit{area and last column generating function} of $R$ we mean the formal sum

\begin{displaymath}
\sum_{P \in R} q^{area\ of\ P} \cdot t^{the\ height\ of\ the\ last\ column\ of\ P}.
\end{displaymath}

\section{Level one polyominoes with cheesy blocks}

Let $E=E(q,t)$ be the area and last column generating function for level one polyominoes with cheesy blocks. Let $E_1=E(q,1)$ and $F_1=\frac{\partial E}{\partial t}(q,1)$. 

Let $U$ be the set of all level one polyominoes with cheesy blocks.

When we build a column-convex polyomino (resp. a cheesy polyomino) from left to right, adding one column at a time, every intermediate figure is a column-convex polyomino (resp. a cheesy polyomino) itself. However, when we build a polyomino with cheesy blocks, this is no longer the case. A ``left factor" of an element of $U$ need not itself be an element of $U$.

We say that a figure $P$ is an \textit{incomplete level one polyomino with cheesy blocks} if $P$ itself is not an element of $U$, but $P$ can be made into an element of $U$ by adding an extra column on the right side of $P$. When we build a level one polyomino with cheesy blocks, then every intermediate object is either a level one polyomino with cheesy blocks or an incomplete level one polyomino with cheesy blocks. Incomplete polyominoes with cheesy blocks can appear when we are building the left ``wing" of some bird. 

Let $V$ be the set of all incomplete level one polyominoes with cheesy blocks. Let

\begin{displaymath}
G(q)=\sum_{P \in V} q^{area \ of \ P}.
\end{displaymath}

For $P \in U \cup V$, we define the \textit{body} of $P$ to be all of $P$, except the rightmost column of $P$.

We write $U_{\alpha}$ for the set of level one polyominoes with cheesy blocks which have only one column. For $P \in U \setminus U_{\alpha}$, we define the \textit{pivot cell} of $P$ to be the lower right neighbour of the lowest cell of the second last column of $P$. As Figure 4 clearly shows, the pivot cell of a polyomino $P \in U \setminus U_{\alpha}$ is not necessarily contained in $P$. Let

\begin{eqnarray*}
U_{\beta} & = & \{P \in U \setminus U_{\alpha}: \mathrm{the \ body \ of \ } P \mathrm{\ lies \ in \ } U\mathrm{, \ the \ last \ column \ of \ } P \\
& & \mathrm{\ has \ no \ hole, \ and \ the \ pivot \ cell \ of \ } P \mathrm{\ is \ contained \ in \ } P \}, \\
U_{\gamma} & = & \{P \in U \setminus U_{\alpha}: \mathrm{the \ body \ of \ } P \mathrm{\ lies \ in \ } U\mathrm{, \ the \ last \ column \ of \ } P \\
& & \mathrm{\ has \ no \ hole, \ and \ the \ pivot \ cell \ of \ } P \mathrm{\ is \ not \ contained \ in \ } P \}, \\
U_{\delta} & = & \{P \in U \setminus U_{\alpha}: \mathrm{the \ body \ of \ } P \mathrm{\ lies \ in \ } U\mathrm{, \ and \ the \ last \ column \ of \ } P \\
& & \mathrm{\ has \ a \ hole\} \quad and} \\
U_{\epsilon} & = & \{P \in U \setminus U_{\alpha}: \mathrm{the \ body \ of \ } P \mathrm{\ lies \ in \ } V \}.
\end{eqnarray*}

\begin{figure}
\begin{center}
\includegraphics[width=118mm]{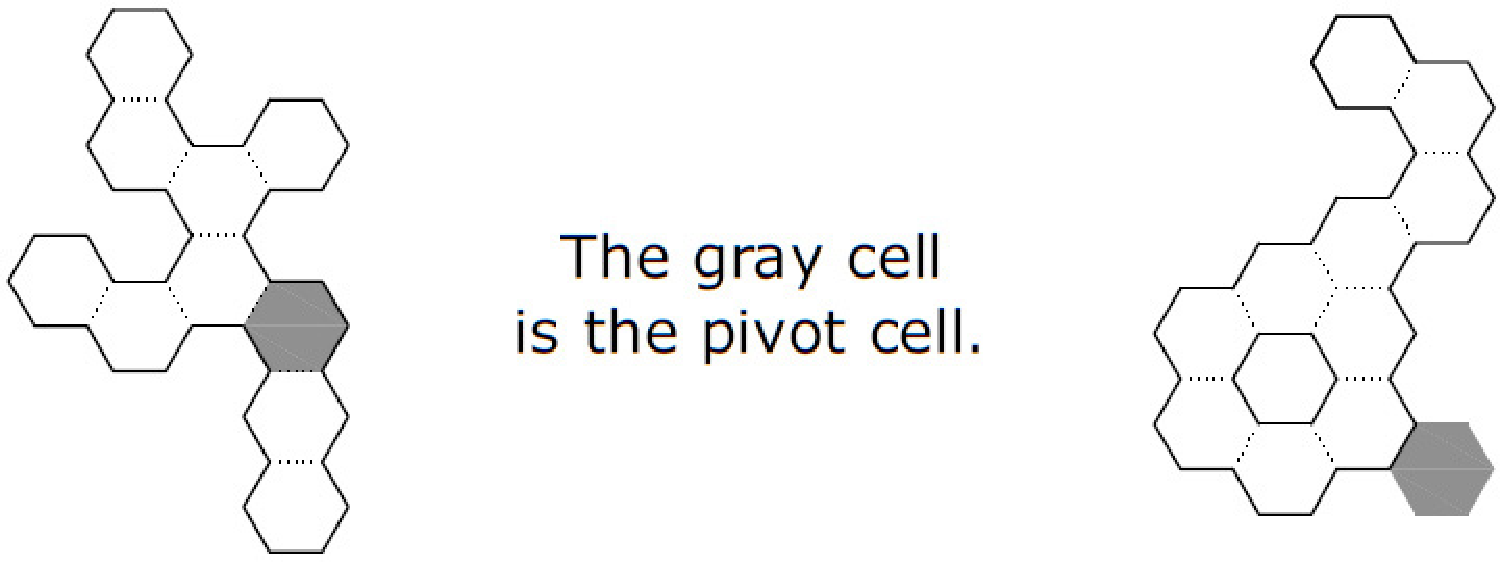}
\caption{The pivot cell.}
\end{center}
\end{figure}

The sets $U_{\alpha}$, $U_{\beta}$, $U_{\gamma}$, $U_{\delta}$ and $U_{\epsilon}$ form a partition of $U$. We write $E_{\alpha}$, $E_{\beta}$, $E_{\gamma}$, $E_{\delta}$ and $E_{\epsilon}$ for the parts of the series $E$ that come from the sets $U_{\alpha}$, $U_{\beta}$, $U_{\gamma}$, $U_{\delta}$ and $U_{\epsilon}$, respectively.

We have

\begin{equation}
E_{\alpha}=qt+(qt)^2+(qt)^3+\ldots =\frac{qt}{1-qt}.
\end{equation}

If a polyomino $P$ lies in $U_{\beta}$, then the last column of $P$ is made up of the pivot cell, of $i \in \{0,\: 1,\: 2,\: 3,\ldots \: \}$ cells lying below the pivot cell, and of $j \in \{0,\: 1,\: 2,\: 3,\ldots \: \}$ cells lying above the pivot cell. See Figure 5. Hence,

\begin{equation}
E_{\beta} = E_1 \cdot qt \cdot \left[ \sum_{i=0}^{\infty} (qt)^i \right] \cdot \left[ \sum_{j=0}^{\infty} (qt)^j \right] = \frac{qt}{(1-qt)^2} \cdot E_1.
\end{equation}

\begin{figure}
\begin{center}
\includegraphics[width=114mm]{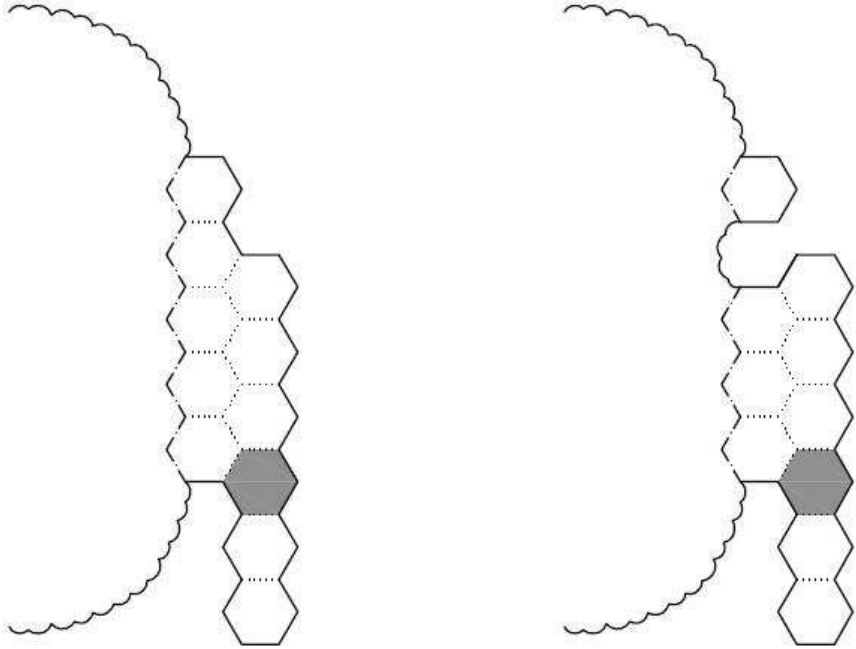}
\caption{The last two columns of two elements of $U_{\beta}$.}
\end{center}
\end{figure}

Consider the following situation. A polyomino $P \in U$ ends with a column $I$. We are creating a new column to the right of $I$, and the result should be an element of $U_{\gamma}$. Then, whether or not the column $I$ has a hole, we can put the lowest cell of the new column in exactly $m$ places, where $m$ is the height of $I$. See Figure 6. Hence

\begin{equation}
E_{\gamma}=\frac{qt}{1-qt} \cdot F_1.
\end{equation}

\begin{figure}
\begin{center}
\includegraphics[width=114mm]{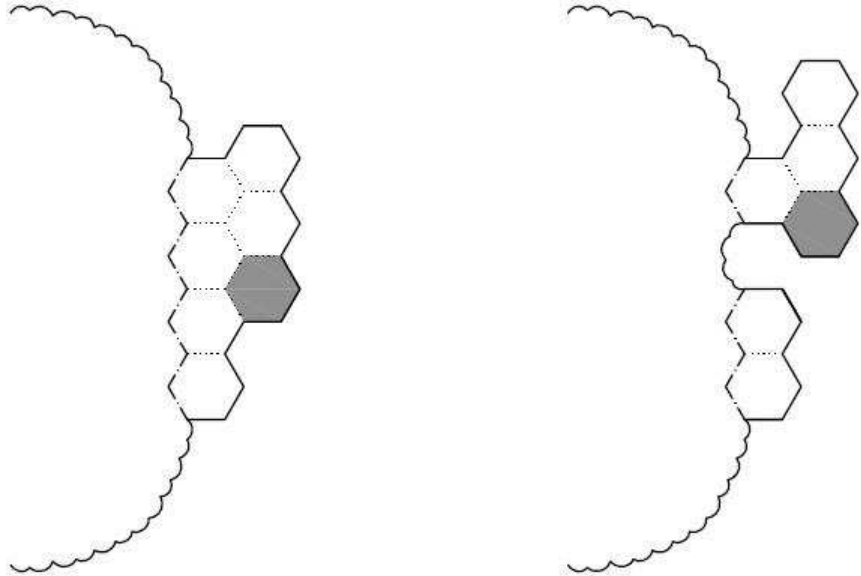}
\caption{The last two columns of two elements of $U_{\gamma}$.}
\end{center}
\end{figure}

Let us proceed to another situation. A polyomino $P \in U$ ends with a column $J$. We are creating a new column to the right of $J$, and the result should be an element of $U_{\delta}$. Then, whether or not the column $J$ has a hole, we can put the hole of the new column in exactly $n-1$ places, where $n$ is the height of $J$. See Figure 7. The new column is made up of $i \in \{1,\: 2,\: 3,\ldots \: \}$ cells lying below the hole, of a hole of height one, and of $j \in \{1,\: 2,\: 3,\ldots \: \}$ cells lying above the hole. Altogether,

\begin{equation}
E_{\delta}=\frac{qt}{1-qt} \cdot t \cdot \frac{qt}{1-qt} \cdot (F_1-E_1) = \frac{q^2 t^3}{(1-qt)^2} \cdot (F_1-E_1).
\end{equation}

\begin{figure}
\begin{center}
\includegraphics[width=114mm]{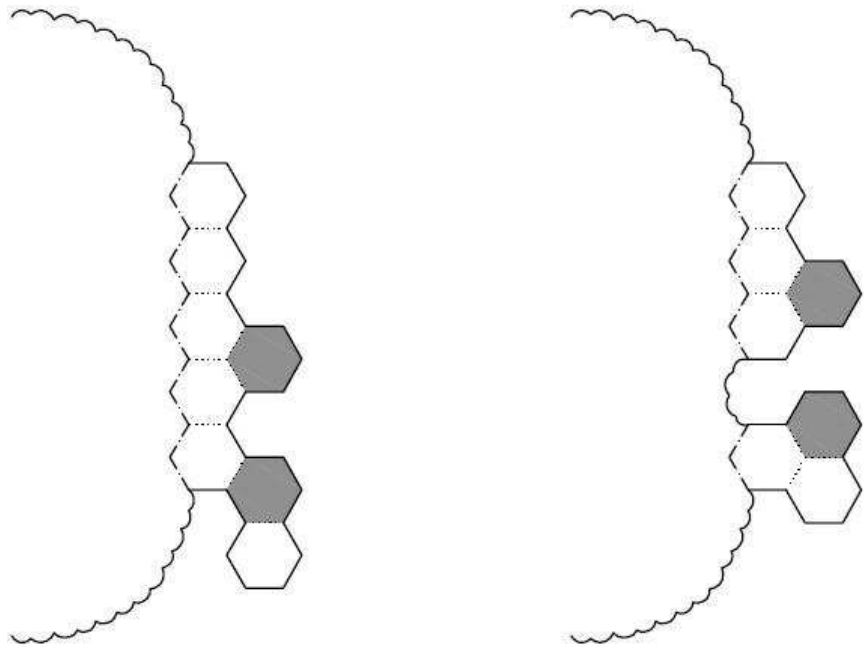}
\caption{The last two columns of two elements of $U_{\delta}$.}
\end{center}
\end{figure}

Now, let $P$ be an element of $U_{\epsilon}$. By the definition of $U_{\epsilon}$, $P$ is a polyomino with cheesy blocks, but the body of $P$ is not a polyomino with cheesy blocks. So, we can decompose $P$ into birds (in one or more ways), but we cannot decompose the body of $P$ into birds. The only possible reason is the following. The body of $P$ ends with some ``problematic" holed columns, while $P$ itself ends with a hole-free column. When this hole-free column is added to the ``problematic" columns, the resulting figure is a bird (with no right ``wing").

Let us translate these remarks into mathematical formulae. In the second last column of $P$ there is a hole, and in the last column of $P$ there are two cells with which the hole is filled. In addition to this two-celled ``cork", the last column contains $i \in \{0,\: 1,\: 2,\ldots \: \}$ cells lying below the ``cork" and 
$j \in \{0,\: 1,\: 2,\ldots \: \}$ cells lying above the ``cork". See Figure 8. Hence

\begin{equation}
E_{\epsilon}=\frac{1}{1-qt} \cdot q^2 t^2 \cdot \frac{1}{1-qt} \cdot G = \frac{q^2 t^2}{(1-qt)^2} \cdot G.
\end{equation}

\begin{figure}
\begin{center}
\includegraphics[width=40mm]{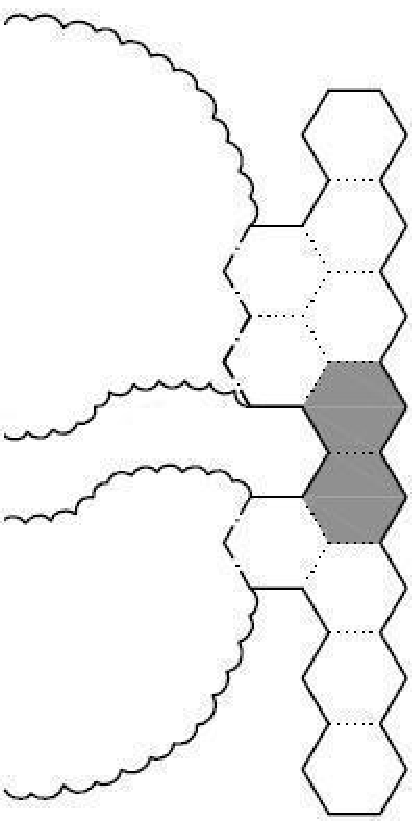}
\caption{The last two columns of an element of $U_{\epsilon}$.}
\end{center}
\end{figure}

Since $E=E_{\alpha}+E_{\beta}+E_{\gamma}+E_{\delta}+E_{\epsilon}$, Eqs. (1)--(5) imply that

\begin{equation}
E=\frac{qt}{1-qt} + \frac{qt}{(1-qt)^2} \cdot E_1 + \frac{qt}{1-qt} \cdot F_1 + \frac{q^2 t^3}{(1-qt)^2} \cdot (F_1-E_1) + \frac{q^2 t^2}{(1-qt)^2} \cdot G.
\end{equation}

Setting $t=1$, from Eq. (6) we get

\begin{equation}
E_1=\frac{q}{1-q} + \frac{q}{(1-q)^2} \cdot E_1 + \frac{q}{1-q} \cdot F_1 + \frac{q^2}{(1-q)^2} \cdot (F_1-E_1) + \frac{q^2}{(1-q)^2} \cdot G.
\end{equation}

Differentiating Eq. (6) with respect to $t$ and then setting $t=1$, we get

\begin{eqnarray}
F_1 & = & \frac{q}{1-q} + \frac{q^2}{(1-q)^2} + \frac{q}{(1-q)^2} \cdot E_1 + \frac{2q^2}{(1-q)^3} \cdot E_1 + \frac{q}{1-q} \cdot F_1 \nonumber \\
& & \mbox{} + \frac{q^2}{(1-q)^2} \cdot F_1 + \frac{3q^2}{(1-q)^2} \cdot (F_1-E_1) + \frac{2q^3}{(1-q)^3} \cdot (F_1-E_1) \nonumber \\
& & \mbox{} + \frac{2q^2}{(1-q)^2} \cdot G + \frac{2q^3}{(1-q)^3} \cdot G.
\end{eqnarray}

Now we turn to incomplete polyominoes with cheesy blocks. Firstly, we see that an incomplete polyomino with cheesy blocks always ends with a holed column.

We write $V_{\alpha}$ for the set of incomplete level one polyominoes with cheesy blocks which have only one column. Let $P \in V \setminus V_{\alpha}$. If the body of $P$ lies in $U$, then the said body is in contact with just one of the two connected components of $P$'s last column. The non-contacting component of the last column is located either wholly above or wholly below the second last column of $P$. We define the \textit{lower pivot cell} of $P \in V \setminus V_{\alpha}$ to be the lower right neighbour of the lowest cell of the second last column of $P$. In addition, we define the \textit{upper pivot cell} of $P \in V \setminus V_{\alpha}$ to be the upper right neighbour of the highest cell of the second last column of $P$. Let

\begin{eqnarray*}
V_{\beta} & = & \{P \in V \setminus V_{\alpha}: \mathrm{the \ body \ of \ } P \mathrm{\ lies \ in \ } U\mathrm{, \ and \ the \ hole \ of \ the \ last} \\
& & \mathrm{column \ of \ } P \mathrm{\ coincides \ either \ with \ the \ lower \ pivot \ cell \ of \ } P \\
& & \mathrm{or \ with \ the \ upper \ pivot \ cell \ of \ } P \} \quad \mathrm{and} \\
V_{\gamma} & = & \{P \in V \setminus V_{\alpha}: \mathrm{the \ body \ of \ } P \mathrm{\ lies \ in \ } U\mathrm{, \ and \ the \ hole \ of \ the \ last} \\
& & \mathrm{column \ of \ } P \mathrm{\ lies \ either \ below \ the \ lower \ pivot \ cell \ of \ } P \\
& & \mathrm{or \ above \ the \ upper \ pivot \ cell \ of \ } P \}.
\end{eqnarray*}

Let us move on to the case when the body of $P \in V \setminus V_{\alpha}$ lies in $V$. Then the second last column of $P$ has two connected components. It is easy to see that each of those two components must be in contact with the last column of $P$. (This does not mean that each of the two connected components of the last column must be in contact with the second last column.) Now, it may or may not happen that one connected component of $P$'s last column is in contact with both connected components of $P$'s second last column. Accordingly, we define the following two sets:

\begin{eqnarray*}
V_{\delta} & = & \{P \in V \setminus V_{\alpha}: \mathrm{the \ body \ of \ } P \mathrm{\ lies \ in \ } V\mathrm{, \ and \ the \ hole \ of \ the \ last} \\
& & \mathrm{column \ of \ } P \mathrm{\ touches \ the \ hole \ of \ the \ second \ last \ column \ of \ } P \} \quad \mathrm{and} \\
V_{\epsilon} & = & \{P \in V \setminus V_{\alpha}: \mathrm{the \ body \ of \ } P \mathrm{\ lies \ in \ } V\mathrm{, \ and \ the \ hole \ of \ the \ last} \\
& & \mathrm{column \ of \ } P \mathrm{\ does \ not \ touch \ the \ hole \ of \ the \ second \ last \ column \ of \ } P \}.
\end{eqnarray*}

The sets $V_{\alpha}$, $V_{\beta}$, $V_{\gamma}$, $V_{\delta}$ and $V_{\epsilon}$ form a partition of $V$. We write $G_{\alpha}$, $G_{\beta}$, $G_{\gamma}$, $G_{\delta}$ and $G_{\epsilon}$ for the parts of the series $G$ that come from the sets $V_{\alpha}$, $V_{\beta}$, $V_{\gamma}$, $V_{\delta}$ and $V_{\epsilon}$, respectively.

The set $V_{\alpha}$ contains every two-part column (with one-celled hole) having $i \in \{1,\: 2,\: 3,\ldots \: \}$ cells below the hole and $j \in \{1,\: 2,\: 3,\ldots \: \}$ cells above the hole. Thus,

\begin{equation}
G_{\alpha} = \frac{q}{1-q} \cdot \frac{q}{1-q} = \frac{q^2}{(1-q)^2}.
\end{equation}

If $P \in V_{\beta}$, then the body of $P$ lies in $U$. The hole of the last column has two possibilities: to coincide with the lower pivot cell of $P$ or to coincide with the upper pivot cell of $P$. Anyhow, the last column is made up of $i \in \{1,\: 2,\: 3,\ldots \: \}$ cells lying below the hole and $j \in \{1,\: 2,\: 3,\ldots \: \}$ cells lying above the hole. See Figure 9. Therefore,

\begin{equation}
G_{\beta} = 2 \cdot \frac{q}{1-q} \cdot \frac{q}{1-q} \cdot E_1 = \frac{2q^2}{(1-q)^2} \cdot E_1.
\end{equation}

\begin{figure}
\begin{center}
\includegraphics[width=114mm]{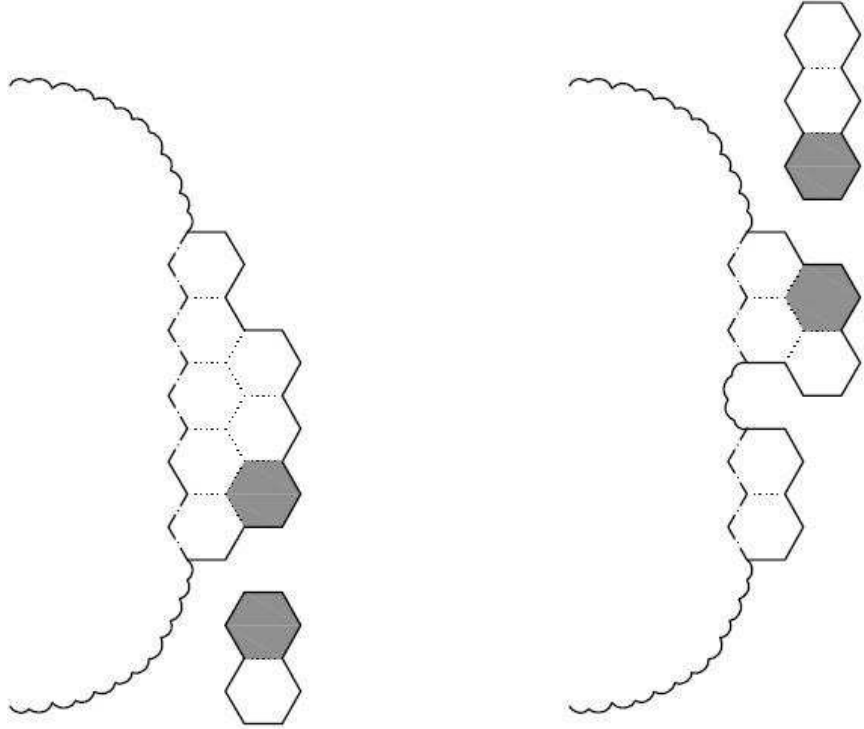}
\caption{The last two columns of two elements of $V_{\beta}$.}
\end{center}
\end{figure}

Now let $P \in V_{\gamma}$. The body of $P$ again lies in $U$. If the hole of the last column lies below the lower pivot cell of $P$, then the last column of $P$ is made up of:

\begin{itemize}
\item $i \in \{1,\: 2,\: 3,\ldots \: \}$ cells lying below the hole,
\item $j \in \{0,\: 1,\: 2,\ldots \: \}$ cells lying above the hole and below the lower pivot cell,
\item the lower pivot cell, and
\item $k \in \{0,\: 1,\: 2,\ldots \: \}$ cells lying above the lower pivot cell.
\end{itemize}

If the hole of the last column lies above the upper pivot cell of $P$, then the last column of $P$ is made up of:

\begin{itemize}
\item $i \in \{1,\: 2,\: 3,\ldots \: \}$ cells lying above the hole,
\item $j \in \{0,\: 1,\: 2,\ldots \: \}$ cells lying above the upper pivot cell and below the hole,
\item the upper pivot cell, and
\item $k \in \{0,\: 1,\: 2,\ldots \: \}$ cells lying below the upper pivot cell.
\end{itemize}

See Figure 10. Altogether,

\begin{equation}
G_{\gamma} = 2 \cdot \frac{q}{1-q} \cdot \frac{1}{1-q} \cdot q \cdot \frac{1}{1-q} \cdot E_1 = \frac{2q^2}{(1-q)^3} \cdot E_1.
\end{equation}

\begin{figure}
\begin{center}
\includegraphics[width=114mm]{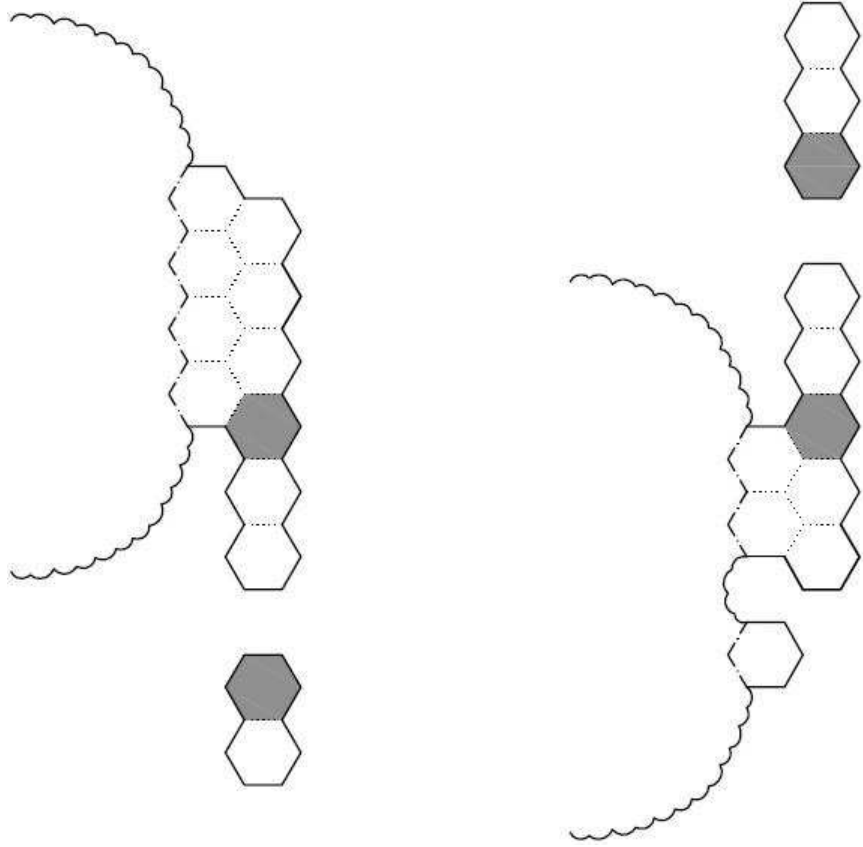}
\caption{The last two columns of two elements of $V_{\gamma}$.}
\end{center}
\end{figure}

If $P \in V_{\delta}$, then the body of $P$ lies in $V$. The second last and last columns of $P$ both have a hole. The hole of the last column is either the lower right neighbour or the upper right neighbour of the hole of the second last column. In the last column, there are $i \in \{1,\: 2,\: 3,\ldots \: \}$ cells below the hole and $j \in \{1,\: 2,\: 3,\ldots \: \}$ cells above the hole. See Figure 11. Hence,

\begin{equation}
G_{\delta} = 2 \cdot \frac{q}{1-q} \cdot \frac{q}{1-q} \cdot G = \frac{2q^2}{(1-q)^2} \cdot G.
\end{equation}

\begin{figure}
\begin{center}
\includegraphics[width=114mm]{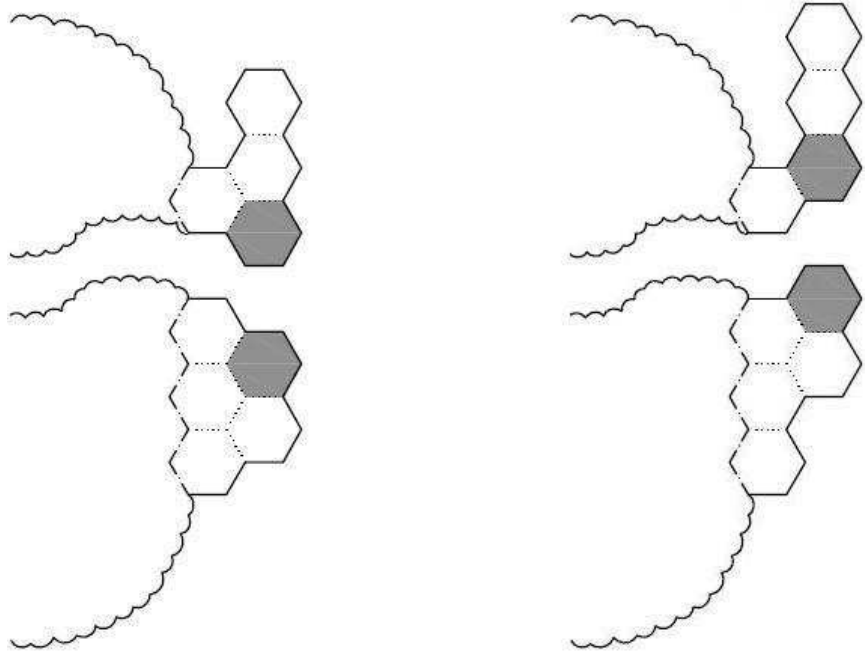}
\caption{The last two columns of two elements of $V_{\delta}$.}
\end{center}
\end{figure}

Let $P \in V_{\epsilon}$. Once again, the second last and last columns of $P$ both have a hole. However, to the right of the hole of the second last column, there are two cells which both belong to $P$.

If this two-celled ``cork" is contained in the upper component of the last column, then the last column is made up of:

\begin{itemize}
\item $i \in \{1,\: 2,\: 3,\ldots \: \}$ cells lying below the hole of the last column,
\item $j \in \{0,\: 1,\: 2,\ldots \: \}$ cells lying above the hole and below the ``cork",
\item the two cells forming the ``cork", and
\item $k \in \{0,\: 1,\: 2,\ldots \: \}$ cells lying above the ``cork".
\end{itemize}

If the said cork is contained in the lower component of the last column, then the last column is made up of:

\begin{itemize}
\item $i \in \{1,\: 2,\: 3,\ldots \: \}$ cells lying above the hole,
\item $j \in \{0,\: 1,\: 2,\ldots \: \}$ cells lying above the cork and below the hole,
\item the two cells forming the cork, and
\item $k \in \{0,\: 1,\: 2,\ldots \: \}$ cells lying below the cork. 
\end{itemize}

See Figure 12. Consequently,

\begin{equation}
G_{\epsilon}=2 \cdot \frac{q}{1-q} \cdot \frac{1}{1-q} \cdot q^2 \cdot \frac{1}{1-q} \cdot G = \frac{2q^3}{(1-q)^3} \cdot G.
\end{equation}

\begin{figure}
\begin{center}
\includegraphics[width=114mm]{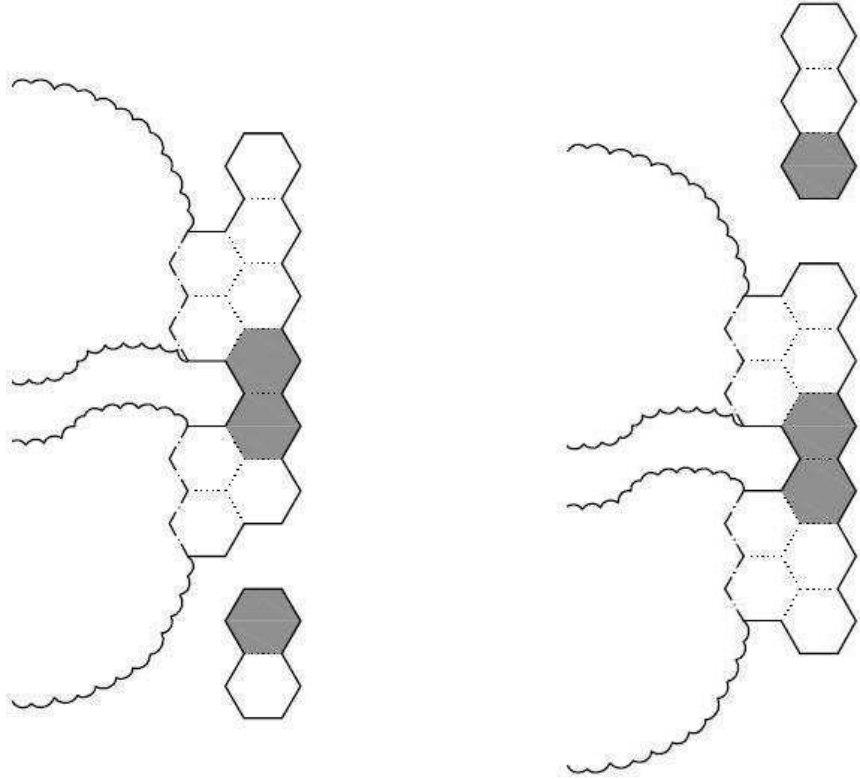}
\caption{The last two columns of two elements of $V_{\epsilon}$.}
\end{center}
\end{figure}

Since $G=G_{\alpha}+G_{\beta}+G_{\gamma}+G_{\delta}+G_{\epsilon}$, Eqs. (9)--(13) imply that

\begin{equation}
G=\frac{q^2}{(1-q)^2} + \frac{2q^2}{(1-q)^2} \cdot E_1 + \frac{2q^2}{(1-q)^3} \cdot E_1 + \frac{2q^2}{(1-q)^2} \cdot G + \frac{2q^3}{(1-q)^3} \cdot G.
\end{equation}

Eqs. (7), (8) and (14) form a system of three linear equations in three un\-knowns, $E_1$, $F_1$ and $G$. Solving this system yields the following result.

\begin{propo} The area generating function for level one polyominoes with cheesy blocks is given by

\begin{displaymath}
E_1=\frac{q(1-6q+11q^2-6q^3+2q^4)}{1-9q+27q^2-32q^3+13q^4-3q^5-q^6} \ .
\end{displaymath}
\end{propo}

The complex roots of the denominator of $E_1$ are\footnote{Those roots have infinitely many digits. Since here we do not have an infinite amount of space, we shall round those roots to six decimal places.} $r_1=-6.109867$, $r_2=0.233117$, $r_3=0.449922-0.087757 \cdot i$, $r_4=0.449922+0.087757 \cdot i$, $r_5=0.988454-1.537589 \cdot i$ and $r_6=0.988454+1.537589 \cdot i$. The root with smallest absolute value is $r_2=0.233117$, and $\frac{1}{r_2}$ is equal to $4.289698$. By decomposing $E_1$ into partial fractions and expanding the partial fractions into Taylor series, we establish the following fact.

\begin{coro} The number of $n$-celled level one polyominoes with cheesy blocks $[q^{n}] E_1$ has the asymptotic behaviour

\begin{displaymath}
[q^{n}] E_1 \sim 0.126651 \cdot 4.289698^n.
\end{displaymath}
\end{coro}

Thus, the growth constant of level one polyominoes with cheesy blocks is $4.289698$. For comparison, the growth constants of level one cheesy polyominoes and column-convex polyominoes are $4.114908$ and $3.863131$, respectively. The increase from $4.114908$ to $4.289698$ is certainly respectable, although not quite so large as the increase from $3.863131$ to $4.114908$.

\section{Level two polyominoes with cheesy blocks}

In this enumeration, if the last column of a polyomino (or of an incomplete polyomino) has two connected components, we often need to record not only the overall height of the last column, but also the height of the last column's upper component and the height of the last column's lower component. Hence, in addition to the ``old" variables $q$ and $t$, we introduce two new variables, $u$ and $v$. As before, the exponent of $q$ is the area and the exponent of $t$ is the overall height of the last column\footnote{Recall what do we mean by the height of a column: in Figure 2, the highlighted column has height $3$, and the next column to the left has height $3$, too.}. The exponent of $u$ is the height of the upper component of the last column, and the exponent of $v$ is the height of the lower component of the last column. 

The four main generating functions in this enumeration are $A=A(q,t)$, $C=C(q,t,u,v)$, $G=G(q,u,v)$ and $J=J(q,u,v)$. Those generating functions are used for the following purposes: 

\begin{itemize}
\item $A$ is a generating function for level two polyominoes with cheesy blocks whose last column either has no hole or has a one-celled hole,
\item $C$ is a generating function for level two polyominoes with cheesy blocks whose last column has a two-celled hole,
\item $G$ is a generating function for incomplete level two polyominoes with cheesy blocks whose last column has a one-celled hole,
\item $J$ is a generating function for incomplete level two polyominoes with cheesy blocks whose last column has a two-celled hole.
\end{itemize}

Let $A_1=A(q,1)$, $B_0=\frac{\partial A}{\partial t}(q,0)$, $B_1=\frac{\partial A}{\partial t}(q,1)$, $C_1=C(q,1,1,1)$, $D_1=\frac{\partial C}{\partial t}(q,1,1,1)$,
$E_0=\frac{\partial C}{\partial u}(q,1,0,1)$, $F_0=\frac{\partial C}{\partial v}(q,1,1,0)$, $G_1=G(q,1,1)$, $H_0=\frac{\partial G}{\partial u}(q,0,1)$, $I_0=\frac{\partial G}{\partial v}(q,1,0)$, $J_1=J(q,1,1)$, $K_0=\frac{\partial J}{\partial u}(q,0,1)$, and $L_0=\frac{\partial J}{\partial v}(q,1,0)$.

In Section 3, where the main generating functions were denoted $E$ and $G$, we established a functional equation Eq. (6) for $E$ and a functional equation Eq. (14) for $G$. Recall that the proof of Eq. (6) breaks into five cases. Namely, we partition the set of all level one polyominoes with cheesy blocks (denoted $U$), into five subsets: $U_{\alpha}$, $U_{\beta}$, $U_{\gamma}$, $U_{\delta}$ and $U_{\epsilon}$. Now, here we have to establish four functional equations, one for each of the generating functions $A$, $C$, $G$ and $J$. The proofs of the functional equations for $A$, $C$, $G$ and $J$ break into nine, two, ten and twelve cases, respectively. It would take quite a lot of space to consider all those cases. Hence, as a kind of compromise, we shall only prove the functional equation for $C$. That will suffice to get a taste of all the four proofs.

Let $S$ be the set of those level two polyominoes with cheesy blocks whose last column either has no hole or has a one-celled hole. Let $T$ be the set of those level two polyominoes with cheesy blocks whose last column has a two-celled hole. As in Section 3, for $P \in T$, we define the \textit{body} of $P$ to be all of $P$, except the rightmost column of $P$. Let

\begin{eqnarray*}
T_{\alpha} & = & \{P \in T: \mathrm{the \ body \ of \ } P \mathrm{\ lies \ in \ } S \}, \quad \mathrm{and} \\
T_{\beta} & = & \{P \in T: \mathrm{the \ body \ of \ } P \mathrm{\ lies \ in \ } T \}.
\end{eqnarray*}

The sets $T_{\alpha}$ and $T_{\beta}$ form a partition of $T$. We write $C_{\alpha}$ and $C_{\beta}$ for the parts of the series $C$ that come from the sets $T_{\alpha}$ and $T_{\beta}$, respectively.

Let $P$ be an element of $S$ and let $c$ be a column with a two-celled hole. Suppose that we want to glue $c$ to $P$ so that $P \cup c$ lies in $T_{\alpha}$, and so that $P$ and $c$ are the body and the last column of $P \cup c$, respectively. In how many ways $P$ and $c$ can be glued together? In principle, the number of ways is (the height of the last column of $P$) minus two. See Figure 13. However, if $P$ ends with a one-celled column, then we can glue $c$ to $P$ in zero ways, and not in minus one ways. Thus, we have

\begin{equation}
C_{\alpha}=\frac{q^2t^4uv}{(1-qtu)(1-qtv)} \cdot (B_1-2A_1+B_0).
\end{equation}

\begin{figure}
\begin{center}
\includegraphics[width=114mm]{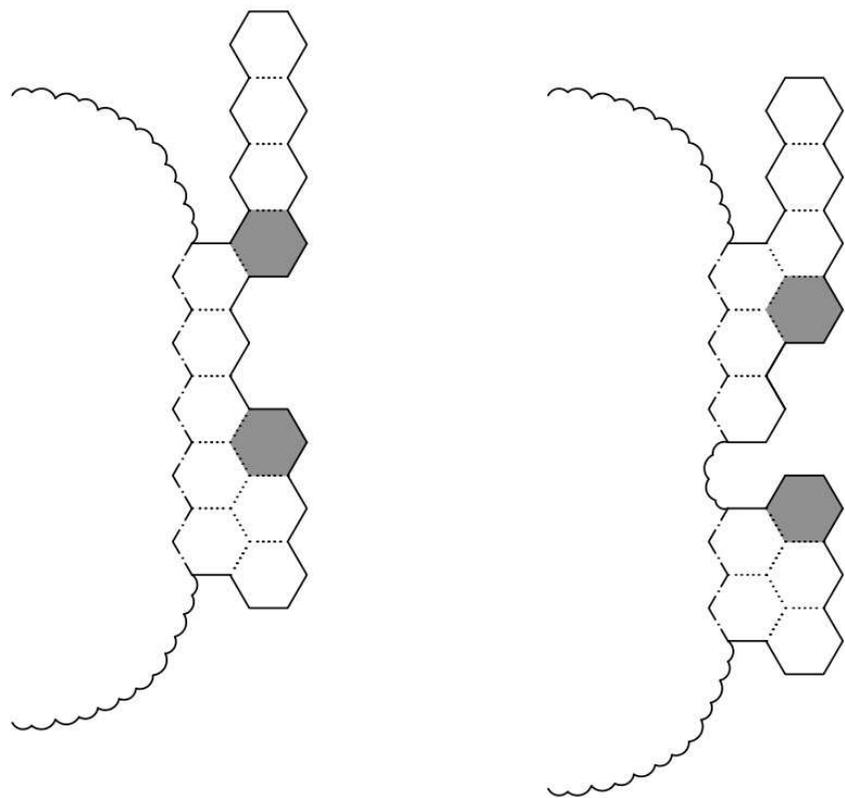}  
\caption{The last two columns of two elements of $T_{\alpha}$.}
\end{center}
\end{figure}

Let $P$ be an element of $T$ and let $c$ be a column with a two-celled hole. Suppose that we want to glue $c$ to $P$ so that $P \cup c$ lies in $T_{\beta}$,
and so that $P$ and $c$ are the body and the last column of $P \cup c$, respectively. In how many ways $P$ and $c$ can be glued together? First, there are (the height of the last column of $P$) minus two ways to satisfy these two necessary conditions:

\begin{itemize}
\item the bottom cell of the upper component of $c$ is either identical with or lies lower than the upper right neighbour of the top cell of the last column of $P$, and
\item the top cell of the lower component of $c$ is either identical with or lies higher than the lower right neighbour of the bottom cell of the last column of $P$. 
\end{itemize}

See Figure 14. 

\begin{figure}
\begin{center}
\includegraphics[width=118mm]{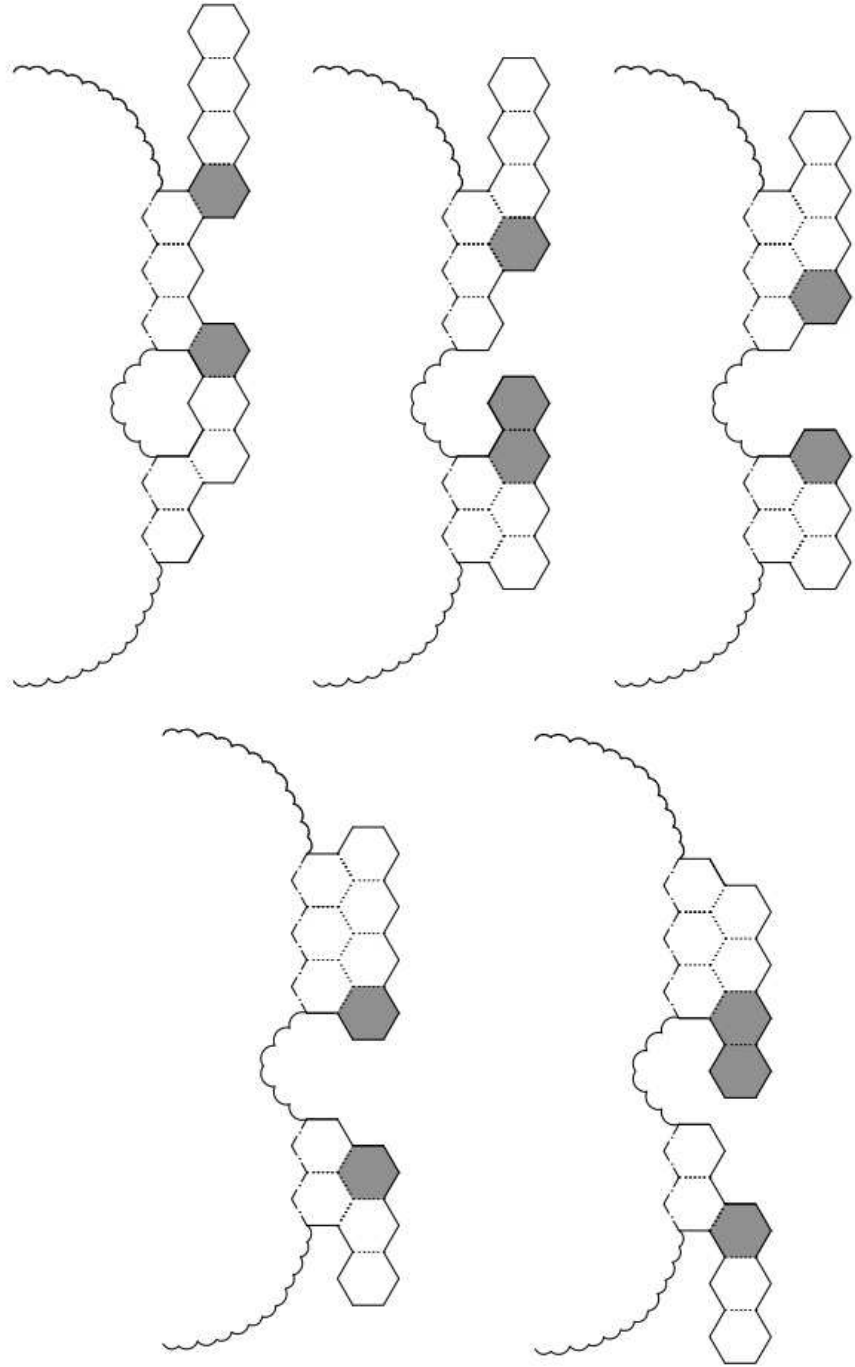}  
\caption{The last two columns of five elements of $T_{\beta}$.}
\end{center}
\end{figure}

So, if there were no special cases, then $C_{\beta}$ would be equal to

\begin{equation}
\frac{q^2t^4uv}{(1-qtu)(1-qtv)} \cdot (D_1-2C_1).
\end{equation}

However, special cases do exist. There are two of them:

\begin{enumerate}
\item The upper component of the last column of $P \in T$ has at least two cells and the lower component of the two-component column $c$ has just one cell.
\item The lower component of the last column of $P \in T$ has at least two cells and the upper component of the two-component column $c$ has just one cell.
\end{enumerate}

In case 1, it is (so to speak) dangerous to glue $c$ to $P$ in such a way that the one-celled lower component of $c$ becomes a common neighbour of the two cells which form the hole of the last column of $P$. This dangerous operation produces an object which is not a polyomino and hence does not lie in $T_{\beta}$.

In case 2, it is dangerous to glue $c$ to $P$ in such a way that the one-celled upper component of $c$ becomes a common neighbour of the two cells which form the hole of the last column of $P$. Again, the dangerous operation produces an object which is not a polyomino and hence does not lie in $T_{\beta}$.

Now, Eq. (16) is actually a generating function for the union of $T_{\beta}$ with the set of objects produced by the two dangerous operations. The generating function for the objects produced by the first dangerous operation is

\begin{equation}
\frac{q^2t^4uv}{1-qtu} \cdot (C_1-E_0).
\end{equation}

The generating function for the objects produced by the second dangerous operation is

\begin{equation}
\frac{q^2t^4uv}{1-qtv} \cdot (C_1-F_0).
\end{equation}

Subtracting Eqs. (17) and (18) from Eq. (16), we obtain

\begin{equation}
C_{\beta}=\frac{q^2t^4uv}{(1-qtu)(1-qtv)} \cdot (D_1-2C_1) - \frac{q^2t^4uv}{1-qtu} \cdot (C_1-E_0) - \frac{q^2t^4uv}{1-qtv} \cdot (C_1-F_0).
\end{equation}

Since $C=C_{\alpha}+C_{\beta}$, Eqs. (15) and (19) imply that

\begin{eqnarray}
C & = & \frac{q^2t^4uv}{(1-qtu)(1-qtv)} \cdot (B_1-2A_1+B_0+D_1-2C_1) \nonumber \\
& & \mbox{}-\frac{q^2t^4uv}{1-qtu} \cdot (C_1-E_0) - \frac{q^2t^4uv}{1-qtv} \cdot (C_1-F_0).
\end{eqnarray}

Setting $t=u=v=1$, from Eq. (20) we get

\begin{equation}
C_1=\frac{q^2}{(1-q)^2} \cdot (B_1-2A_1+B_0+D_1-2C_1) - \frac{q^2}{1-q} \cdot (2C_1-E_0-F_0).
\end{equation}

Differentiating Eq. (20) with respect to $t$ and then setting $t=u=v=1$, we get

\begin{eqnarray}
D_1 & = & \left[\frac{4q^2}{(1-q)^2} + \frac{2q^3}{(1-q)^3}\right] \cdot (B_1-2A_1+B_0+D_1-2C_1) \nonumber \\
& & \mbox{} - \left[\frac{4q^2}{1-q}+\frac{q^3}{(1-q)^2}\right] \cdot (2C_1-E_0-F_0).
\end{eqnarray}

Next, we differentiate Eq. (20) with respect to $u$ and then set $t=1$, $u=0$ and $v=1$. The result is

\begin{equation}
E_0=\frac{q^2}{1-q} \cdot (B_1-2A_1+B_0+D_1-3C_1+F_0) -q^2 \cdot (C_1-E_0).
\end{equation}

Also, we differentiate Eq. (20) with respect to $v$ and then set $t=1$, $u=1$ and $v=0$. The result is

\begin{equation}
F_0=\frac{q^2}{1-q} \cdot (B_1-2A_1+B_0+D_1-3C_1+E_0) -q^2 \cdot (C_1-F_0).
\end{equation}

Eqs. (21)--(24) make part of a system of altogether $13$ linear equations in $13$ unknowns: $A_1$, $B_0$, $B_1$, $C_1$, $D_1$, $E_0$, $F_0$, $G_1$, $H_0$, $I_0$, $J_1$, $K_0$ and $L_0$. The other $9$ equations of that linear system are obtained from the functional equations for $A$, $G$ and $J$. The computer algebra system \textit{Maple} quickly solved the linear system and then summed the generating functions $A_1$ and $C_1$. The result can be seen in the following proposition.

\begin{propo}
The area generating function for level two polyominoes with cheesy blocks is given by

\begin{displaymath}
M=\frac{N}{O},
\end{displaymath}

\noindent where

\begin{eqnarray*}
N & = & q \cdot (1-13q+70q^2-202q^3+336q^4-317q^5+143q^6+18q^7-84q^8 \\
& & \mbox{}+11q^9+227q^{10}-375q^{11}+267q^{12}-165q^{13}+134q^{14}-21q^{15}+4q^{16} \\
& & \mbox{}-124q^{17}+98q^{18}-12q^{19}+28q^{20}-16q^{21})
\end{eqnarray*}

\noindent and

\begin{eqnarray*}
O & = & 1-16q+107q^2-391q^3+850q^4-1108q^5+797q^6-169q^7-266q^8 \\
& & \mbox{}+317q^9+159q^{10}-913q^{11}+1081q^{12}-672q^{13}+446q^{14}-268q^{15} \\
& & \mbox{}+7q^{16}-158q^{17}+404q^{18}-222q^{19}+42q^{20}-70q^{21}+34q^{22}.
\end{eqnarray*}
\end{propo}

\begin{coro}
The number of $n$-celled level two polyominoes with cheesy blocks has the asymptotic behaviour

\begin{displaymath}
[q^{n}] M \sim 0.102214 \cdot 4.462811^n.
\end{displaymath}
\end{coro}

So, the growth constant of level two polyominoes with cheesy blocks is $4.462811$.

\section{Level three polyominoes with cheesy blocks}

In this section, we skip everything but the final results.

\begin{propo}
The area generating function for level three polyominoes with cheesy blocks is given by

\begin{displaymath}
P=\frac{Q}{R},
\end{displaymath}

\noindent where

\begin{eqnarray*}
Q & = & q \cdot (1-24q+264q^2-1766q^3+8033q^4-26297q^5+63860q^6 \\
& & \mbox{}-116445q^7+157849q^8-148533q^9+61825q^{10}+99443q^{11} \\
& & \mbox{}-308464q^{12}+519182q^{13}-655900q^{14}+618461q^{15}-344081q^{16} \\
& & \mbox{}-101610q^{17}+519331q^{18}-707969q^{19}+601249q^{20}-284943q^{21} \\
& & \mbox{}-68043q^{22}+297023q^{23}-346370q^{24}+265550q^{25}-140577q^{26} \\
& & \mbox{}+31503q^{27}+64681q^{28}-166424q^{29}+234520q^{30}-218182q^{31} \\
& & \mbox{}+130432q^{32}-29144q^{33}-33391q^{34}+38482q^{35}-12237q^{36}-2050q^{37} \\
& & \mbox{}-6144q^{38}+18593q^{39}-21514q^{40}+11634q^{41}+3351q^{42}-13907q^{43} \\
& & \mbox{}+12096q^{44}+2302q^{45}-8825q^{46}+570q^{47}+4681q^{48}-1695q^{49} \\
& & \mbox{}-1519q^{50}+1290q^{51}+64q^{52}-224q^{53}+44q^{54}-12q^{55}) 
\end{eqnarray*}

\noindent and

\begin{eqnarray*}
R & = & 1-27q+334q^2-2515q^3+12906q^4-47836q^5+132248q^6 \\
& & \mbox{}-276956q^7+438796q^8-508406q^9+365771q^{10}+36865q^{11} \\
& & \mbox{}-648120q^{12}+1344653q^{13}-1932847q^{14}+2126787q^{15}-1632701q^{16} \\
& & \mbox{}+408884q^{17}+1117382q^{18}-2223607q^{19}+2392085q^{20}-1636807q^{21} \\
& & \mbox{}+418146q^{22}+665251q^{23}-1211688q^{24}+1191386q^{25}-838060q^{26} \\
& & \mbox{}+416174q^{27}-41907q^{28}-323733q^{29}+664097q^{30}-810808q^{31} \\
& & \mbox{}+657803q^{32}-319442q^{33}+14159q^{34}+120746q^{35}-95202q^{36} \\
& & \mbox{}+22341q^{37}-7930q^{38}+47294q^{39}-74720q^{40}+62640q^{41}-19120q^{42} \\
& & \mbox{}-28394q^{43}+46822q^{44}-21864q^{45}-18416q^{46}+20930q^{47}+6617q^{48} \\
& & \mbox{}-14093q^{49}+982q^{50}+5867q^{51}-2682q^{52}-642q^{53}+608q^{54} \\ 
& & \mbox{}-88q^{55}+12q^{56}.
\end{eqnarray*}
\end{propo}

\begin{coro}
The number of $n$-celled level three polyominoes with cheesy blocks has the asymptotic behaviour

\begin{displaymath}
[q^{n}] P \sim 0.090504 \cdot 4.538766^n.
\end{displaymath}
\end{coro}

\section{Taylor expansions and the limit value of the growth constants}

To see how many polyominoes of a given type have $1,\: 2,\: 3,\ldots$ cells, we expanded the area generating functions into Taylor series. The results are shown in Table~1. In Table 2, we display the growth constant of column-convex polyominoes, together with all the growth constants which we computed in this paper and in \cite{semi}.

\begin{table}
\begin{center}
\begin{tabular}{|r||r|r|r|r|}\hline
& \multicolumn{1}{c|}{Column-} & \multicolumn{1}{c|}{Level 1} & \multicolumn{1}{c|}{Level 2} & \multicolumn{1}{c|}{Level 3} \\
& \multicolumn{1}{c|}{convex} & \multicolumn{1}{c|}{polyominoes} & \multicolumn{1}{c|}{polyominoes} & \multicolumn{1}{c|}{polyominoes} \\
& \multicolumn{1}{c|}{polyo-} & \multicolumn{1}{c|}{with cheesy} & \multicolumn{1}{c|}{with cheesy} & \multicolumn{1}{c|}{with cheesy} \\
\multicolumn{1}{|c||}{Area} & \multicolumn{1}{c|}{minoes} & \multicolumn{1}{c|}{blocks} & \multicolumn{1}{c|}{blocks} & \multicolumn{1}{c|}{blocks} \\ \hline \hline
1 & 1 & 1 & 1 & 1 \\ \hline
2 & 3 & 3 & 3 & 3 \\ \hline
3 & 11 & 11 & 11 & 11 \\ \hline
4 & 42 & 44 & 44 & 44 \\ \hline
5 & 162 & 184 & 186 & 186 \\ \hline
6 & 626 & 784 & 810 & 812 \\ \hline
7 & 2419 & 3363 & 3582 & 3614 \\ \hline
8 & 9346 & 14451 & 15952 & 16259 \\ \hline
9 & 36106 & 62097 & 71242 & 73558 \\ \hline
10 & 139483 & 266716 & 318441 & 333683 \\ \hline
11 & 538841 & 1145074 & 1423411 & 1515454 \\ \hline
12 & 2081612 & 4914448 & 6360809 & 6885303 \\ \hline
\end{tabular}
\caption{Here is how many polyominoes of a given type have $1,\: 2,\ldots,\: 12$ cells.}
\end{center}
\end{table}

\begin{table}
\begin{center}
\begin{tabular}{|c||c|c|} \hline
& & Polyominoes \\
& Cheesy & with cheesy \\
Level & polyominoes & blocks \\ \hline \hline
0 & 3.863131 & 3.863131 \\ \hline
1 & 4.114908 & 4.289698 \\ \hline
2 & 4.231836 & 4.462811 \\ \hline
3 & 4.288631 & 4.538766 \\ \hline
\end{tabular}
\caption{The growth constants. By level $0$ cheesy polyominoes, and so too by level $0$ polyominoes with cheesy blocks, we mean the usual column-convex polyominoes.}
\end{center}
\end{table}

Now, it is natural to ask the question: to what value do the growth constants tend as level tends to infinity? Our database is too small for giving a precise answer. Anyway, we shall permit ourselves to make a vague estimate. In the case of polyominoes with cheesy blocks, computing the first differences of the growth constants, we obtain the numbers $4.290-3.863=0.427$, $4.463-4.290=0.173$, and $4.539-4.463=0.076\: $. The sequence $0.427,\: 0.173,\: 0.076$ is a little similar to a geometric sequence with common ratio $\frac{2}{5}$. So, the limit value of the growth constants of polyominoes with cheesy blocks might be about 
$4.463 + \frac{5}{3} \cdot 0.076 = 4.590\: $. In a similar way, in \cite{semi}, the limit value of the growth constants of cheesy polyominoes was estimated to be about $4.346\: $.

\end{document}